\documentclass[12pt]{article}
\usepackage{mathrsfs}
\usepackage{amsfonts}
\usepackage{bbm}
\usepackage{amssymb}
\usepackage{amsmath,amssymb}
\newcommand{\al}{\alpha}

\newcommand{\eset}{\emptyset}

\newcommand{\G}{\Gamma}

\newcommand{\mf}{\mathfrak}
\newcommand{\nin}{\noindent}

\newcommand{\se}{\subset}
\newcommand{\seq}{\subseteq}

\newcommand{\vs}{\vspace*}

\def\Lla{\Longleftarrow}
\def\Lra{\Longrightarrow}
\def \nin {\noindent}

\def \Lemma #1 {\vs{3mm}\nin {\bf Lemma #1} \it}
\def \Prop #1 {\vs{3mm}\nin {\bf Proposition #1} \it}
\def \Th #1 {\vs{3mm}\nin {\bf Theorem #1} \it}
\def \Cor #1 {\vs{3mm}\nin {\bf Corollary #1} \it}
\def \Ex #1 {\vs{3mm}\nin {\bf Example #1} \it}
\def \Proof {\vs{3mm}\nin {\bf Proof. }}
\def \part #1 {\hfil\break\hglue 12pt {\rm (#1)~}}

\def \qed {~\vrule height6pt width 6pt depth 0pt}
\def\fs{\footnotesize}

\title{
\bf\LARGE On zero divisors and prime elements of bounded semirings\thanks{This research is supported by the National Natural
Science Foundation of China (Grant No. 11271250).}}
\author{{Tongsuo Wu$^a$\thanks{Corresponding author,\, tswu@sjtu.edu.cn (T. Wu)}, Yuanlin Li$^b$\thanks{yli@brocku.ca (Y. Li)},
Dancheng Lu$^c$\thanks{ludancheng@suda.edu.cn (D. Lu)}}\\
  {\small $^a$Department of Mathematics, Shanghai Jiaotong University}\\
  {\small Shanghai 200240, P. R. China}\\
  {\small $^b$Department of Mathematics, Brock University, St. Catharines, On., Canada L2S 3A1}\\
   {\small $^c$Department of Mathematics, Suzhou University, Suzhou 215006, P.R. China}\\
  \date{}
}
\begin{document}
\baselineskip=16pt
\maketitle

\begin{center}
\begin{minipage}{12cm}

 \vs{3mm}\nin{\small\bf Abstract.} {\fs  A bounded semiring $A$ is a poset together with a compatible  semiring structure. A bounded semiring has the important property of being locally semimodular as a poset, thus the chain complex of the poset $A$ is shellable. In this paper,  properties of zero divisors and prime elements of a bounded semiring are studied. In particular, it is proved that under some mild assumption the set $Z(A)$ of nonzero zero divisors of $A$ is $A\setminus \{0,1\}$, each  prime element of $A$ is a maximal element.
 For a  bounded semiring $A$ with $Z(A)=A\setminus \{0,1\}$, it is proved that $A$ has finitely many maximal elements if ACC holds either for  elements of $A$ or for principal annihilating ideals of $A$. As applications of prime elements,  the structure of a bounded semiring $A$ is completely determined by the structure of integral bounded semirings if  either $|Z(A)|=1$ or  $|Z(A)|=2$ and $Z(A)^2\not=0$.  Applications to the ideal structure of commutative rings are considered.
}

\vs{3mm}\nin {\small Key Words:} {\small bounded semiring, zero divisor, prime element, small $Z(A)$, ideal structure of rings}

\end{minipage}
\end{center}

\vs{4mm}\nin{\bf 1. Introduction and preliminaries }

\vs{3mm}\nin   Throughout this paper, all semigroups $S$ and all rings $R$ are assumed to be commutative with zero element $0_S$ ( i.e., $0S=0$) and  with identity $1_R$ respectively. For a semigroup $S$, let $Z(S)$ be the set of nonzero zero divisors and $S^*$  the set of (nonzero) elements of $S$. For a ring $R$,  let $U(R)$ be the set of invertible elements of $R$, $N(R)$ the nil radical  and $J(R)$ the Jacobson radical of $R$.

 A {\it commutative semiring} is a set $A$ which contains at least two elements $0,1$ and which is equipped with two binary operations, $+$ and $\cdot$, called addition and multiplication respectively, such that the following conditions hold:

1. $(A, +,0)$ is a commutative monoid with zero element $0$.

2. $(A, \cdot,1)$ is a commutative monoid with identity element 1.

3. Multiplication distributes over addition.

4. $0$ annihilates $A$ with respect to multiplication, i.e., $0¡¤a = 0, \,\forall a \in A$.

If there is no zero-divisor in a semiring $A$, then $A$ is called an {\it integral semiring}. Certainly, each ring is a semiring. Other important examples of semirings include the set $\mathbb I(R)$ of ideals of a commutative ring $R$,  the set $\mathbb N$ of nonnegative integers, and the real segment $[0,1]$ whose addition is $max$ operation.  Note that both $\mathbb N$ and $[0,1]$ are integral semirings. Based on a semiring structure, an excellent and rather general framework for constraint satisfaction and optimization was developed in \cite{BMR,BMRSVF}.

Next we introduce a new notion which will be the central topic of this paper.

\vs{3mm}\nin{\bf Definition 1.1.} A {\it partially-ordered semiring} (abbreviated as a {\it po-semiring}) is a commutative semiring $(A, +, \cdot,0,1 )$, together with a compatible partial order $\le$, i.e., a partial order $\le$ on the underlying set $A$ that is compatible with the semiring operations in the sense that it satisfies the following conditions:

(1) $x\le y$ implies $x+z\le y+z$, and

(2) $0\le x$ and $y\le z$  imply that $xy\le xz$
for all $x,y,z$ in $A$.

\nin If $A$ satisfies the following additional condition, then $A$ is called a {\it bounded semiring}:

(3) The partial set $(A,\le, 0,1)$ is bounded, i.e., $1$ is the largest element and $0$ is the least element of $A$.

\vs{3mm}We  remark  that  condition  (3)  is  so  strong  that  it  forces  a  bounded semiring  $A$  to  be  a  dioid, where a semiring is called a {\it dioid} if  its addition is idempotent ($ a+a=a,\,\forall a\in A$).
Furthermore, the above defined partial order $\le$ for a bounded semiring $A$ is identical with the new partial order $\le_1 $ defined by the following
$$\text{$a\le_1 b$ if and only if $a+b=b$}.$$
In other words, a bounded semiring $(A,+,0,1)$ is a bounded join-semilattice.
Clearly, any bounded, distributive lattice is a commutative dioid under join and meet, where $a\le b$ iff $a\wedge b=a$. Each bounded, distributive lattice is certainly also a bounded semiring under the Definition $1.1$.

We also remark that a bounded semiring $A$ is a {\it bounded, locally semimodular poset, thus the chain complex of the poset $A$
is shellable}, see \cite{Bjorner} for details.

An element $p$ of a  bounded semiring is called {\it prime}, if $p\not=1$ and $xy\le p$ implies either $x\le p$ or $y\le p$. An element $x$ is called {\it minimal}, if $x\not=0$ and $0<y\le x$ implies $x=y$ . An element $\mathfrak{m}$ is called {\it maximal}, if $\mathfrak{m}\not=1$ and $\mathfrak{m}\le x<1$ implies $\mathfrak{m}=x$.

An {\it ideal} $I$ of a bounded semiring $A$ is an additive sub-semigroup containing the zero element $0_A$ such that $IA\seq I$.  For an element $u$ of a bounded semiring $A$, set $$<u>\,=\{x\in A\,|\,x\le u\}.$$Clearly,  $<u>$ is an ideal of $A$, called the {\it lower principal ideal generated by $u$}. The principal ideal $uA$ is certainly another ideal generated by $u$, and clearly $uA\seq\, <u>$. An ideal $I$ of $A$ is called {\it hereditary}, if $<u>\,\seq I$ holds for all $u$ in $I$. For any
element $u$ of $A$, both $<u>$ and $ann_A(u)$ are hereditary ideals of $A$, where $ann_Au=\{x\in A\,|\,xu=0\}$.

Our prototype of bounded semiring is the bounded semiring $\mathbb I(R)$ of a commutative ring $R$, which consists of all ideals of $R$. The multiplication is the ideal multiplication, and the addition is the addition of subsets. The partial order is the usual inclusion. The bounded semiring $\mathbb I(R)$ also has the property that for each element $v\not=1$, there exists a maximal element $\mathfrak{m}$ such that $v\le \mathfrak{m}.$ $\mathbb I(R)$ also has the completeness property, where a bounded semiring $A$ is   called to have {\it completeness property} if for each nonempty subset $S$ of $A$, there is the notion of sum $\sum_{x\in S}x$ in $A$.  Note that $\mathbb I(R)$ is also a bounded lattice under $I\wedge J=I\cap J$ and $I\vee J=I+J$, but it is not necessarily a distributive lattice for a general ring $R$. Denote  $\{0,1\}=\mathbb I(F)$ for any field $F$, and note that it is the smallest bounded semiring, i.e., it can be embedded into any bounded semiring.

Throughout the paper, we always assume that one of the following three additional conditions on a bounded semiring $A$ holds:

\vs{2mm}$(C_1)$: {\it For each non-nilpotent element $u$ of $A$, there is a nonzero idempotent $w$ such that $w\le u$ and $w$ has an orthogonal idempotent complement in $A$, i.e., there exists an idempotent element $v$ in $A$ such that $w+v=1,wv=0$.}
\vs{2mm}

\vs{2mm}$(C_2)$: {\it For each  nonzero idempotent element $u$ of $A$, there is a nonzero idempotent $w$ such that $w\le u$ and $w$ has an orthogonal idempotent complement in $A$.}

\vs{2mm}$(C_3)$: {\it  Each idempotent minimal element
of $A$ has an orthogonal idempotent complement in $A$.}
\vs{2mm}

Clearly,  $(C_1)\Lra (C_2)\Lra(C_3)$ hold  for any bounded semiring $A$. In sections 2 and 3, examples will be given to show that the reverse implications do not hold. Note that if $Z(A)\not=\eset$, then each minimal element of $A$ is a zero divisor (see Lemma 3.4(2)).

\vs{3mm}\nin{\bf Proposition 1.2.} {\it Let $R$ be a commutative ring and let let $A=\mathbb{I}(R)$. Then

$(1)$ $A$ satisfies condition $(C_3)$.}

(2) {\it If $R$ is a noetherian exchange ring, then $A$ satisfies condition $(C_2)$.}

(3)   {\it If $R$ is a noetherian exchange ring and $J(R)=N(R)$, then $A$ satisfies condition $(C_1)$. In particular, it holds for any  artinian ring $R$.}

(4) {\it If $A$ satisfies condition $(C_1)$, then $J(R)=N(R)$.}

\vs{3mm}\Proof (1) It follows from Brauer's Lemma, see \cite[10.22]{Lam}.

 (2) Let $I$ be a nontrivial idempotent ideal of $R$. Then $I$ is finitely generated, and thus $I\not\seq J(R)$ by Nakayama Lemma. Then the exchange property ensures the existence of a nonzero idempotent $e$ in $I$. Thus $(Re)^2=Re$ and it has the orthogonal idempotent complement $R(1-e)$. Hence condition $(C_2)$ holds for $\mathbb I(R)$.

(3) An artinian ring $R$ is a noetherian exchange ring and each non-nilpotent ideal of $R$ contains a nonzero idempotent.

(4) For any $x\not\in N(R)$, $Rx$ is not nilpotent. By condition $(C_1)$, there exists a nonzero idempotent $e$ in $Rx$. Thus $x\not\in J(R)$ and it shows $J(R)=N(R)$. \quad\qed

\vs{3mm}We remark that exchange rings include von Neumann regular rings, artinian rings, semilocal rings such that idempotents lift modulo their Jacobson radical. For further information on exchange rings, see \cite{LW} and the listed  references.

In this paper we investigate properties of zero divisors and prime elements of a bounded semiring.   Section 2 deals with the following questions: (1) Under what conditions can it occur that $Z(A)=A\setminus \{0,1\}$? (2) When does a bounded semiring $A$ have only finitely many  maximal elements? (3) When is every prime element also maximal?  It is proved that $Z(A)=A\setminus \{0,1\}$, each  prime element of $A$ is a maximal element, if one of the following conditions is satisfied: (i) Condition $(C_2)$ holds in $A$, and DCC holds for elements of $A$. (ii) Condition $(C_1)$ holds in $A$, and there exists no infinite set of orthogonal idempotents in $A$. For a  bounded semiring $A$ with $Z(A)=A\setminus \{0,1\}$, it is also proved that $A$ has finitely many maximal elements if ACC holds either for  elements of $A$ or for principal annihilating ideals of $A$.

In Section 3, we study the structure of a bounded semiring $A$ with $1\le |Z(A)|\le 2$. It is shown that the structure of a bounded semiring $A$ is completely determined by the structures of integral bounded semirings, if  either $|Z(A)|=1$ or $|Z(A)|=2$ and $Z(A)^2\not=0$ (Theorems 3.2 and 3.3). When $A$ is taken to be $\mathbb I(R)$ for some commutative ring $R$, applications to the ideal structure of a ring are provided.

\vs{4mm}\nin{\bf 2. Chain conditions on bounded semirings }

\vs{3mm}\nin In this section, we study properties of elements of $A\setminus \{0,1\}$. We begin with an easy proposition which will be used repeatedly.

\vs{3mm}\nin {\bf Proposition 2.1.} {\it Let $A$ be a bounded semiring. Then }

(1) {\it Each maximal element of $A$ is prime.}

(2) {\it For a maximal element $\mathfrak{m}$ of $A$ and any element $b$ of $A$, if $\mathfrak{m}b=0$, then either $\mathfrak{m}^2=\mathfrak{m}$ or $b^2=0$.}

(3) {\it If $e_i+f_i=1$ where $e_if_i=0,\,e_i^2=e_i,\,f_i^2=f_i$, then $e_1>e_2$ implies $f_1<f_2$.}

\vs{3mm}\Proof (1) Let $\mathfrak{m}$ be a maximal element of $A$. Assume that $ab\le \mathfrak{m}$ and $a\not\le \mathfrak{m}$. Then $a+\mathfrak{m}> \mathfrak{m}$, so $a+\mathfrak{m}=1$. Then $b=ab+\mathfrak{m}b\le\mathfrak{m}+\mathfrak{m}=\mathfrak{m},$ so $\mathfrak{m}$ is a prime element of $A$.

(2) For a maximal element $\mathfrak{m}$ of $A$ and any $b\in A$, either $\mathfrak{m}+b=\mathfrak{m}$ or $\mathfrak{m}+b=1$. If $\mathfrak{m}+b=\mathfrak{m}$, then $b^2=0$. If $\mathfrak{m}+b=1$, then $\mathfrak{m}=\mathfrak{m}^2.$

(3) If $e_1\ge e_2$, then $e_1+f_2\ge e_2+f_2=1 $, and hence $f_1\le f_2f_1\le f_2.$ If $f_1=f_2$, then $e_1=e_1(e_2+f_2)=e_1e_2\le e_2$. This completes the verification. \quad\qed

\vs{3mm}\nin {\bf Theorem 2.2.} {\it For a bounded semiring A, $Z(A)=A\setminus \{0, 1\}$ if one of the following conditions holds:}

(1) {\it Condition $(C_1)$ holds and there exists no infinite set of orthogonal idempotents in
$A$.}

(2) {\it DCC holds for elements of $A$, and condition $(C_2)$ holds.}

{\it In each case, for any element $c$ in $A\setminus \{0, 1\}$, either $c$ is nilpotent or there exist a positive integer $n$ and a nontrivial idempotent $e$ such that $c^n=c^ne$, where $e$ has an orthogonal idempotent complement in $A$. }

\vs{3mm}\Proof (1) Suppose that condition $(C_1)$ holds in $A$. Suppose to the contrary that  there exists an element $a$ in  $A\setminus \{0,1\}$ such that $a\not\in Z(A)$. Then $a$ is not nilpotent. Hence
by condition $(C_1)$, there exist nonzero idempotents $e_1,f_1$ such that $e_1\le a, e_1f_1=0,e_1+f_1=1.$ Clearly, $af_1$ is not nilpotent. ( Note that if $af_1$ is nilpotent and $a$ is idempotent, then $af_1=0$ and thus $a=ae_1\le e_1$, so $a=e_1$. This will be used to obtain Theorem 2.3.)

 Since $af_1$ is not nilpotent by assumption on $a$, by condition $(C_1)$, there exist nonzero orthogonal idempotents $e_2,f_2$ such that $e_2+f_2=1,e_2\le af_1$. Then $e_1e_2=0$, so  $e_1<e_1+e_2$ and $(e_1+e_2)+(f_1f_2)=1$.  Note that $f_1f_2$ is a nonzero idempotent element and is orthogonal to the idempotent $(e_1+e_2)$. Clearly, $a(f_1f_2)$ is not nilpotent. (If $a(f_1f_2)$ is nilpotent and $a$ is idempotent, then $af_1f_2=0$, implying $a=a(e_1+e_2)\le e_1+e_2$, so $a=e_1+e_2$.)

Since $a(f_1f_2)$ is not nilpotent by assumption on $a$,  there exist nonzero orthogonal idempotents $e_3,f_3$ such that $e_3+f_3=1,e_3\le a(f_1f_2)$. Then
we have an orthogonal idempotent decomposition $(e_1+e_2+e_3)+(f_1f_2f_3)=1$, and thus $e_1<e_1+e_2<e_1+e_2+e_3$, where $e_ie_j=\delta_{ij}e_i$. (If $a$ is idempotent and $a(f_1f_2f_3)=0$, then
$a=e_1+e_2+e_3$. Clearly, $a(f_1f_2f_3)$ is not nilpotent under the assumption on $a$.)

Continuing this process,  we finally obtain  an infinite set $\{e_1,e_2,\cdots\}$ of orthogonal idempotents in $A$, contradicting the assumption on $A$. Note that $e_1<e_1+e_2<e_1+e_2+e_3<\cdots$ and $(e_1+\cdots+e_i)+(f_1\cdots f_i)=1$ imply that $f_1>f_1f_2>\cdots$ by Proposition 2.1(3). (Note also that $e_1+\cdots +e_i\le a$. Note that if assume $a(f_1f_2\cdots f_i)=0$ for some $i$, then $a=e_1+e_2+\cdots+e_i$.)

(2) Assume that condition $(C_2)$ holds for $A$ and  DCC holds for elements of $A$.  Suppose to the contrary that in  $A\setminus \{0,1\}$ there exists an element $b$ such that $b\not\in Z(A)$. Then $b$ is not nilpotent and there exists a positive integer $m$ such that $b^m$ is idempotent. Let $a=b^m$ in the proof of (1). By repeating the discussions in the proof in part (1), we obtain an infinite descending chain of idempotent elements $f_1>f_1f_2>\cdots$, giving a contradiction. \quad \qed

\vs{3mm} Recall that a nonzero idempotent of a semiring is called a {\it primitive idempotent} if it cannot be written as a sum of two orthogonal nontrivial idempotents. By the proof of Theorem 2.2, we have the following improved result for idempotent elements.

\vs{3mm}\nin{\bf Theorem 2.3.} {\it  For a  bounded semiring $A$, if condition $(C_2)$ holds in $A$ and $A$ contains no infinite set of orthogonal idempotents, then each idempotent element of $A$ has an orthogonal idempotent complement.   In particular, each nontrivial idempotent is a zero divisor of $A$. Furthermore, each nonzero idempotent of $A$ is a finite sum of orthogonal primitive idempotents.}

\Proof  The result follows from the proof of (1) in Theorem 2.2, if we start with a nontrivial idempotent $a$. \quad \qed

\vs{3mm}\vs{3mm}We remark that there exists no infinite set of orthogonal idempotents in a bounded semiring  $A$ provided that one of the following conditions holds: (1) ACC holds for idempotent elements of $A$. (2) DCC holds for idempotent elements of $A$, and condition $(C_2)$ also holds. The latter conclusion follows from Theorem 2.3 and Proposition 2.1(3). If ACC (respectively, DCC) holds for elements of $A$, then for each element $a\not=0,1$, clearly there is a maximal (respectively, minimal) element $x$ of $A$ such that $x\ge a$ (respectively, $x\le a$).

 Let $A=\mathbb{I}(R)$ for some commutative ring $R$. In view of Theorems 2.2 and 2.3, we have the following applications to commutative rings.

\vs{3mm}\nin{\bf Corollary 2.4.} {\it Assume that a commutative ring $R$ satisfies one of the following conditions:}

(1) {\it Each non-nilpotent ideal of $R$ contains a nonzero idempotent element, and $R$ contains no infinite set of orthogonal idempotent elements.}

(2) {\it $R$ is artinian. }

 {\nin\it Then for each nontrivial ideal $I$ of $R$, there exists a nontrivial idempotent element $e$ such that $I\seq Re$. Furthermore,  each nonzero idempotent ideal of $R$ has the form $\sum_{i=1}^rRe_i$, where $e_1,\cdots,e_r$ are orthogonal primitive idempotents.}

 \vs{3mm} An ideal $I$ of a semiring $A$ is called an {\it annihilating} ideal if $I=ann_A(S)$ for some nonempty subset $S$ of $A$.   Call an ideal $I$ of $A$ a {\it principal annihilating ideal} if $I=ann_A(u)$ for some element $u$ of $A$.

\vs{3mm}\nin{\bf Corollary 2.5.} {\it For a commutative noetherian  exchange ring $R$, each nontrivial idempotent ideal of $R$ is an annihilating ideal. Furthermore, each nonzero idempotent ideal of $R$ has the form $\sum_{i=1}^rRe_i$, where $e_1,\cdots,e_r$ are orthogonal primitive nonzero idempotents. }

\Proof By applying Proposition 1.2(1) and Theorem 2.3 to $\mathbb{I}(R)$, we obtain the result.\quad\qed

\vs{3mm}Note that any artinian ring is a noetherian exchange ring. Thus Corollary 2.5 holds in particular for artinian rings.

It is well-known that in a commutative artinian ring $R$, each prime ideal is a maximal ideal of $R$, i.e., each prime element of the bounded semiring $\mathbb{I}(R)$ is a maximal element  if DCC holds for elements of $\mathbb{I}(R)$. Also, if DCC holds for elements of $\mathbb{I}(R)$, then ACC also holds. The following example shows that the above mentioned results are not true for a general bounded semiring. It also shows that the additional condition $(C_2)$ in Theorem 2.3 is needed, and that $Z(A)=A\setminus \{0,1\}$ does not imply condition $(C_2)$.

\vs{3mm}\nin{\bf Example  2.6} {\it There exists an infinite bounded semiring $A$ such that DCC holds but ACC fails for elements of $A$, $Z(A)=A\setminus \{0,1\}$ and $A$ has infinitely many prime elements none of which is a maximal element.}

 Let $A=\{0,1,a,b_1,b_2,\cdots\}$ be a countable set with $4\le |A|\le \aleph_0$, and define a partial order $\le$ by $0<a<b_1<b_2<\cdots<1$ on $A$. Define an addition by $x+y=\text{max\{$x,y$\}}, \,\forall x,y\in A$. Define a commutative multiplication by $$0x=0,\,1x=x\,(\forall x\in A),\,ab_i=0,\,a^2=0,\,b_ib_j=b_{min\{i,j\}}.$$ Then it is routine to check that $A$ is a bounded semiring. Note that $a$ and $b_i$ are prime elements of $A$ for all $i$, and  condition $(C_2)$ does not hold for $A$ although condition $(C_3)$ does hold.  The  elements of $A$ satisfies DCC, but they do not satisfy ACC  if $|A|$ is infinite. Note also that $Z(A)=A\setminus \{0,1\}$, and there exists no infinite set of orthogonal idempotents.

\vs{3mm}{\it On the other hand, for a noetherian ring $R$ which is not artinian, ACC holds for elements of the bounded semiring $\mathbb{I}(R)$ but DCC fails.}

\vs{3mm} Despite of the above example,  we are able to show that under some suitable condition the set of the maximal elements is the same as the set of the prime elements of $A$.

\vs{3mm}\nin{\bf Theorem 2.7} {\it For a bounded semiring $A$, each prime element of $A$  is a maximal element if one of the following condition holds: }

(1) {\it Condition $(C_2)$ holds in $A$, and DCC holds for elements of $A$.}

(2) {\it Condition $(C_1)$ holds in $A$, and there exists no infinite set of orthogonal idempotents in $A$.}

\Proof (1) Suppose to the contrary that there exists a prime element $q$ such that $q<p$ for some $p<1$. Then $p$ is not nilpotent. By the DCC assumption, there exists a positive integer $r$ such that $p^r$ is nonzero and idempotent. By Theorem 2.3, $p^r$ is annihilated by a nonzero idempotent, say $f$.

By condition $(C_2)$, there exist orthogonal nonzero idempotents $f_1,g_1$ such that $f_1\le f$, $g_1+f_1=1$. Then $p^r\cdot f_1=0$. Since $q<p$ and $q$ is prime, we have $f_1\le q$  and hence  $g_1\not\le q$. Clearly, $p^rg_1=p^r\not\le q$ and it is idempotent.

By condition $(C_2)$, there exist nonzero orthogonal idempotents $t_2,s_2$ such that $t_2\le p^rg_1$, $t_2+s_2=1$. Then we have $f_1t_2=0$. This together with $t_2+s_2=1$ and $f_1+g_1=1$ implies $f_1+t_2+g_1s_2=1,$ where $f_1,t_2,g_1s_2$ are mutually orthogonal idempotents. If $s_2\not\le q$, then $t_2\le q$. In this case, let $f_2=t_2,\,g_2=g_1s_2$. If $s_2\le q$, then $t_2\not\le q$ and $g_1s_2\le q$. In this case, let $f_2=g_1s_2,g_2=t_2$. In either case, we have an orthogonal idempotent decomposition $f_1+f_2+g_2=1$, where $f_1\le q,\,f_2\le q$ and $g_2\not\le q$.

The next step is to consider the idempotent $p^rg_1g_2$. Clearly $p^rg_1g_2\not\le q$ and in particular, $p^rg_1g_2\not=0$. By condition $(C_2)$, there exist orthogonal idempotents $t_3,s_3$ such that $t_3+s_3=1$ and $t_3\le p^rg_1g_2$. Clearly, $t_3f_i=0$ for $i=1,2$. This together with  $(f_1+f_2)+g_2=1$ implies $f_1+f_2+t_3+s_3g_2=1$. Since $f_1+f_2+t_3\le p$, $s_3g_2\not=0$. If $s_3\not\le q$, then $t_3\le q$. In this case, let $f_3=t_3,\,g_3=s_3g_2$. If $s_3\le q$, then $t_3\not\le q$ and $g_2s_3\le q$. In this case let $f_3=g_2s_3,\,g_3=t_3$. In either case, we have an orthogonal idempotent decomposition $(f_1+f_2+f_3)+g_3=1$, where $f_i\le q$ ($\forall i=1,2,3$) and $g_3\not\le q$.

Continuing this process, we have got an infinite set of mutually orthogonal idempotents $\{f_1,f_2,\cdots\}$. Since $f_1<f_1+f_2<f_1+f_2+f_3<\cdots$, by Proposition 2.1(3) we have an infinite descending chain of idempotents $g_1>g_2>\cdots, $ contradicting the DCC assumption.

(2) Suppose to the contrary that there exists a prime element $q$ such that $q<p$ for some $p<1$. Then $p$ is not nilpotent. By Theorem 2.2, there exists a positive integer $r$ such that $p^r$ is annihilated by a nonzero idempotent, say, $f$. The rest of the proof is almost the same as in (1). Note that $p^r(g_1g_2\cdots g_s)\not\le q$ implies that $p^r(g_1g_2\cdots g_n)$ is not nilpotent and thus condition $(C_1)$ applies. \quad\qed

\vs{3mm}\nin{\bf Corollary 2.8.}   {\it For a commutative ring $R$, each prime ideal of $R$ is a maximal ideal if one of the following condition is satisfied:}

(1) {\it There exists no infinite set of orthogonal idempotents in $R$ and each non-nilpotent ideal contains a nonzero idempotent.}

(2) {\it $R$ is an exchange ring, $J(R)=N(R)$ and there exists no infinite set of orthogonal idempotents in $R$. }

(3) {\it $R$ is an artinian ring.}

\vs{3mm}   The proof of the following theorem is a typical example of use of ideas in proving the Chinese Remainder Theorem (see, e.g., $[4]$).

\vs{3mm}\nin{\bf Theorem  2.9.}  {\it For a  bounded semiring $A$, assume $Z(A)=A\setminus \{0,1\}$. If ACC holds either  for principal annihilating ideals of $A$ or for  elements of $A$, then $A$ has finitely many maximal elements. }

\vs{3mm}\Proof Suppose that $Z(A)=A\setminus \{0,1\}$. Suppose to the contrary that $A$ has infinitely many maximal elements and let $\mathfrak{m}_i$ ($i\in N^{+}$) be distinct maximal elements of $A$. Then $\mathfrak{m}_1+\mathfrak{m}_3=1$ and $\mathfrak{m}_2+\mathfrak{m}_3=1$ and hence $1=\mathfrak{m}_3(\mathfrak{m}_1+\mathfrak{m}_2+\mathfrak{m}_3)+\mathfrak{m}_1\mathfrak{m}_2= \mathfrak{m}_3+\mathfrak{m}_1\mathfrak{m}_2$. In a similar way,  we obtain $\mathfrak{m}_1\mathfrak{m}_2\cdots \mathfrak{m}_m+\mathfrak{m}_{m+1}=1$. Now consider the following ascending chain of principal annihilating ideals of the semiring $A$
$$ann_A(\mathfrak{m}_1)\seq\cdots \seq ann_A(\mathfrak{m}_1\cdots \mathfrak{m}_m)\seq ann_A(\mathfrak{m}_1\cdots \mathfrak{m}_m\mathfrak{m}_{m+1})\seq \cdots.$$

If ACC holds for principal annihilating ideals of $A$, then there exists an integer $m$ such that $ann_A(\mathfrak{m}_1\cdots \mathfrak{m}_m)= ann_A(\mathfrak{m}_1\cdots \mathfrak{m}_m\mathfrak{m}_{m+1})$. By assumption, there is a nonzero element $x\in ann_A(\mathfrak{m}_{m+1})$. But then $x=x\mathfrak{m}_{m+1}+x(\mathfrak{m}_1\cdots \mathfrak{m}_m)=0$, a contradiction.

If ACC holds for elements of $A$, we claim that there exists no strict ascending chain $ann_A(\mathfrak{m}_1)<ann_A(\mathfrak{m}_1\mathfrak{m}_2)<\cdots $, and the result thus follows. In fact, if this were not the case, then for each $m$ there would exist an element $y_{m+1}$ such that $$y_{m+1}\in ann_A(\mathfrak{m}_1\cdots \mathfrak{m}_m\mathfrak{m}_{m+1})\setminus ann_A(\mathfrak{m}_1\cdots \mathfrak{m}_m).$$
Then we would have obtained an infinite ascending  chain of elements of $A$: $$y_2<y_2+y_3<y_2+y_3+y_4<\cdots.$$ This completes the proof. \qed

\vs{3mm} Combining Theorems 2.2, 2.7 and 2.9, we have the following results.

\vs{3mm}\nin{\bf Corollary 2.10. }
 {\it Let $A$ be a bounded semiring satisfying condition $(C_2)$. If both ACC and DCC holds for elements of $A$, then $A$ has only a finite number of prime elements. These primes are precisely the maximal elements of $A$. Moreover, $Z(A)=A\setminus \{0,1\}$.}

 \vs{3mm}\nin{\bf Corollary 2.11. } {\it Let $A$ be a bounded semiring satisfying condition $(C_1)$. If ACC holds for elements of $A$, then $A$ has only a finite number of maximal elements.}

 \vs{3mm} We remark that if $(C_1)$ holds in a bounded semiring $A$ and $A$ has a unique maximal element $\mf{m}$,   then $\mathfrak{m}$ is nilpotent.

 \vs{3mm} Applying Theorem 2.9 to $\mathbb{I}(R)$ for a commutative ring $R$, we have the following.

\vs{3mm}\nin{\bf Corollary 2.12.}  {\it Let $R$ be a commutative noetherian ring.}

(1) (\cite[[Proposition 1.7]{BR}) {\it If each nontrivial ideal is an annihilating ideal, then $R$ is a semilocal ring. }

(2) {\it If each non-nilpotent ideal of $R$ contains an idempotent, then $R$ is a semilocal ring.}

(3) {\it Any Noetherian exchange ring $R$ with $J(R)=N(R)$ is a semilocal ring.}

\vs{3mm}Recall that an ideal $I$ of  a bounded semiring $A$ is called {\it hereditary}, if $<u>\,\seq I$ holds for all $u$ in $I$.  Motivated by Theorem 2.9 we have the following.

\vs{3mm}\nin{\bf Proposition 2.13} {\it Let $A$ be a bounded semiring. Then}

 (1) {\it  ACC holds for elements of $A$ if and only if  ACC holds for hereditary ideals of $A$.}

 (2) {\it If  DCC holds for hereditary ideals of $A$, then DCC also holds for elements of $A$. The converse holds, if  $A$ has the completeness property and any hereditary ideal $I$ is closed under taking infinite sums. }

 \Proof
  $\Lla$ of $(1)$ and $(2)$: For any $u\in A$,  $<u>$ is clearly a hereditary ideal of the bounded semiring $A$. For elements $u,v$ of $A$, $u\le v$ if and only if
$<u>\,\seq\, <v>$, and $u< v$ if and only if $<u>\,\subset\, <v>$. This implies the sufficiency part of the proposition.

$\Lra $: (1) If ACC does not hold for hereditary ideals of  the bounded semiring $A$, then there exists a strict ascending chain of hereditary ideals of $A$, say, $X_1\se X_2\se \cdots$. Then for any $n\ge 2$, take an element $u_n$ of $R$ such that $u_n\in X_n\setminus X_{n-1}$. Clearly, there is an infinite ascending chain $u_1< u_1+u_2< u_1+u_2+u_3< \cdots,$ and hence ACC does not hold in $A$.

(2) If DCC does not hold for hereditary ideals of $A$, then there exists a strict descending chain of hereditary ideals of $A$, say, $Y_1\supset Y_2\supset \cdots$. Then for any $n\ge 1$, take a $v_n$ such that $v_n\in Y_n\setminus Y_{n+1}$. Clearly, there is an infinite descending chain of elements of $A$, namely $\sum_{i\ge 1} v_i> \sum_{i\ge 2} v_i>\cdots,$  where $\sum_{i\ge m} v_i\in Y_m$. Then DCC does not hold for elements of $A$. This completes proof.\quad \qed

\vs{3mm}\nin{\bf Corollary 2.14.} (1) {\it If ACC (DCC, respectively)  holds for ideals of a bounded semiring $A$, then ACC (DCC, respectively) also holds for elements of $A$.}

(2) {\it If ACC holds for elements of $A$, then ACC holds for principal annihilating ideals of $A$.}

(3) {\it For a bounded semiring $A$, assume that  $A$ has the completeness property  and any hereditary ideal $I$ is closed under taking infinite sums. If further DCC holds for elements of $A$, then DCC holds for principal annihilating ideals of $A$. }

\vs{3mm}\nin{\bf Corollary 2.15.}  {\it Let $R$ be a commutative ring and denote $A=\mathbb I(R)$. Then ACC (DCC, respectively) holds for elements of $A$ if and only if  ACC (DCC, respectively) holds for hereditary ideals of $A$. }

\vs{3mm}The proof of the following result is routine and is omitted here.

\vs{3mm}\nin{\bf Proposition 2.16} {\it For a bounded semiring $A$ and an element $p$ of $A$, $p$ is a prime element of $A$ if and only if $<p>$ is a prime ideal of $A$.}

\vs{4mm}\nin{\bf 3. The structure of a bounded semiring $A$ with small $|Z(A)|$}

\vs{3mm}In this section, we study the structure of a bounded semiring $A$ with small $Z(A)$.  We first have the following lemma.

\vs{3mm}\nin{\bf Lemma 3.1.}  {\it Let $A$ be a bounded semiring.}

(1) {\it If $Z(A)=\{c\}$, then $c^2=0$, $c$ is the least nonzero element of $A$ and is a prime element of $A$.}

(2) {\it If $Z(A)=\{c,u\}$, then exactly one of the following holds:}

\hspace{0.5cm} (2.1) {\it Both c and u are minimal elements of $A$, $c<x,\,u<x,\, \forall x\in A\setminus\{0,c,u\}$, $u^2=u$ and $ c^2=c$. Furthermore, both $c$ and $u$ are prime elements of $A$.}

\hspace{0.5cm} (2.2) {\it $c$ is the least nonzero element of $A$, and  $c^2=0$. $u$ is a prime element of $A$, and $u<p$ for each prime element $p\not=c,u$. }

\Proof (1) If $Z(A)=\{c\}$, then clearly $c^2=0$. Now if $xy\le c$ for some $x,y\in A$, then $x(xy^2)=0$. If further $x\not\le c$, then $xy^2=0$ and hence $y^2=0$. Then $y\le c$. This shows that $c$ is a prime element of $A$. Now for any $x\in A\setminus \{0,c\}$, clearly $c=xc$ and hence $c<x$. Thus $c$ is the least nonzero element of $A$.

(2) Let $Z(A)=\{c,u\}$. First, assume that $c$ and $u$ are incomparable. Then  $c+u\not\in \{c,u\}$ and hence $c^2=c,u^2=u$. It follows that $xc=c,\,xu=u,\,\forall x\not=0,c,u$, and hence $c<x,u<x$ hold for each $x\not=0,c,u$. It implies that both $c$ and $u$ are minimal elements of $A$. To verify that $c$ is a prime element, assume $xy\le c$ and $y\not\le c$. If  $y=u$, then $x\le c$. If further $y\not=u$, then $0=(xy)u=xu$ and so $x\le c$. Thus $c$ is a prime element of $A$. By symmetry, $u$ is also a prime element of $A$.

Next we assume $c< u$. Then $c^2=0$. For any $x\in A\setminus \{0,c,u\}$, clearly $c=xc$ and hence $c<x$. Thus $c$ is the least nonzero element of $A$. For any $x\not\le u,y\not\le u$, $xy\not\le u$ since otherwise, $0=(xy)c=x(yc)=c$, a contradiction. Thus $u$ is a prime element of $A$. Finally, for any prime element $p\not=c,u$, if $u^2\in \{0,c\}$, then $u^3=0$ and hence $u< p$. If $u^2=u$, then $pu\not=c$ and hence $u=pu< p$.\quad \qed

\vs{3mm}\nin{\bf Theorem 3.2.} {\it $A$ is a bounded semiring with $|Z(A)|=1$ if and only if  there exist an integral bounded semiring $A_1$ and an element $c\not\in A_1$ such that $A\cong \{c\}\cup A_1$, where the partial order of $ \{c\}\cup A_1$ is extended from $A_1$ by adding $0<c<a$ $(\forall a\in A_1^*)$, and the commutative addition and multiplication in $\{c\}\cup A_1$ are extended respectively from that of the bounded semiring $A_1$, and the following conditions are fulfilled:
$$0+c=c,\, c+c=c,\, c^2=0,\,c+y=y,\, cy=c \,(\forall y\in A_1^*).$$}

\Proof  $\Lra$: Let $A$ be a bounded semiring and denote $A_1=A\setminus Z(A)$. Then $A_1$ is an integral sub-bounded semiring of $A$, and  the rest of the results follow from Lemma 3.1(1).

$\Lla$: It is routine to check that $\{c\}\cup A_1$ is a bounded semiring under the assumptions. Clearly $|Z(A)|=1.$ \quad\qed

\vs{3mm}\nin{\bf Theorem 3.3.} {\it  $A$  is a bounded semiring such that $|Z(A)|=2$ and $Z(A)^2\not=0$ if and only if one of the following holds:}

 (1)  {\it There exist an integral bounded semiring $A_1$ and elements $c,u\not\in A_1$ such that $A\cong \{c,u\}\cup A_1$, where $A_1$ has a least nonzero element $a_0$, the partial order of
 $\{c,u\}\cup A_1$ is extended from $A_1$ by adding $0<c<x$ and $0<u<x$ $(\forall x\in A_1^*)$, and the commutative addition and multiplication in $\{c\}\cup A_1$ are extended respectively from that of the bounded semiring $A_1$, and the following conditions are fulfilled:
 $$c+u=a_0,\,0+x=x,\,x+x=x,\, x+y=y\,(\forall x\in \{c,u\},\, \forall y\in A_1^*),$$
$$u^2=u,\,c^2=c,\, cu=0=0x,\,xy=x\, (\forall x\in\{c,u\},\,\forall y\in A_1^*).$$
}

 (2) {\it There exist an integral bounded semiring $A_1$ and elements $c,u\not\in A_1$ such that $A\cong \{c,u\}\cup A_1$, where the partial order of $ \{c,u\}\cup A_1$ is extended from $A_1$ by adding $0<c<u<a$ $(\forall a\in A_1^*)$, and the commutative addition and multiplication in $\{c\}\cup A_1$ are extended respectively from that of the bounded semiring $A_1$, and the following conditions are fulfilled:
 $$x+y=max\{x,y\},\, \forall x\in \{c,u\},\,\forall y\in \{c,u\}\cup A_1,$$
$$c^2=0,\,cu=0=0x,\, u^2\in \{c,u\},\,xy=x\, (\forall x\in\{c,u\},\,\forall y\in A_1^*).$$
}

\Proof $\Lra$: Assume that $A$  is a bounded semiring such that $|Z(A)|=2$. Assume further $Z(A)=\{c,u\}$, and set $A_1=A\setminus \{c,u\}$.

(1) If $c$ and $u$ are incomparable, then $c+u\in A_1^*$ and   $Z(A)^2\not=0$ by Lemma 3.1(2). Assume $c+u=a_0$. By the proof of Lemma $3.1(2)$, we have $a_0y=(c+u)y=c+u=a_0, \forall y\in A_1^*$. Thus $a_0$ is the least nonzero element of $A_1$. Since both $c$ and $u$ are prime elements of $A$ by Lemma 3.1(2), $A_1$ is an integral sub-bounded semiring of $A$.  The other statements also follow directly from Lemma 3.1(2).

(2) In the following, assume that $c$ and $u$ are comparable, assume further $c<u$. Then $c^2=0$, and $c$ is the least nonzero element of $A$. Thus $ac=c$ for any nonzero element $a$ in $A_1^*$. If $Z(A)^2\not=0$, then $u^2\not=0$. This implies $au=u, \forall a\in A_1^*$. Then $u<a,\,\forall a\in A_1^*$. Since $u$ is a prime element of $A$, $u\not\in A_1^*A_1^*$. Clearly, $c\not\in A_1^*A_1^*$. Thus $A_1$ is an integral sub-bounded semiring of $A$. This completes the necessary part of the proof.

$\Lla$: It is not hard to check that $\{c,u\}\cup A_1$ is a bounded semiring under either assumption. Clearly, $|Z(A)|=2$ and $Z(A)^2\not=0$ hold in either case. Note that $(C_3)$ holds in  Case (2).\quad\qed

\vs{3mm}  We remark that  Theorem 3.3 (1) gives a complete characterization of bounded semirings $A$ with $|Z(A)|=2$, in which there exist no nilpotent elements. In particular, we have the following.

\vs{3mm}\nin {\bf Corollary 3.4} (1) {\it A bounded semiring $A$ satisfies condition $(C_3)$, $|Z(A)|=2$ and $Z(A)^2\not=0$ if and only if either $A\cong \{0,1\}\times \{0,1\}$, or $A$ is the bounded semiring constructed in Theorem 3.3(2).}

(2){\it A bounded semiring $A$ satisfies condition $(C_3)$, $|Z(A)|=2$ and there exist no nilpotent elements in $A$ if and only if $A\cong \{0,1\}\times \{0,1\}$.}

\vs{3mm}Note that the bounded semiring $\mathbb{I}(R)$ always satisfies condition $(C_3)$. Applying Corollary 3.4 to the bounded semiring $\mathbb{I}(R)$,  we have the following Corollary 3.5 which is essentially the same with \cite[Theorem 3]{AA2} by \cite[Theorem 2.1]{McL}. To find definition and some results on the annihilating ideal graph $\mathbb{AG}(R)$, the reader is referred to $[1,2,5]$. Note that the graph $\mathbb{AG}(R)$ defined in \cite{BR} is exactly the zero divisor graph of the multiplicative semigroup $\mathbb{I}(R)$. Thus all known results on graph properties for the semigroup graph $\G(S)$ hold for the annihilating ideal graph $\mathbb{AG}(R)$. One can refer to $[3,8,12]$ for some further fundamental properties of the zero-divisor graph $\G(R)$ of a semigroup $R$ or a ring $R$.  We remark that in $[14]$ we classified all bipartite  zero-divisor graphs $\G(A)$ (respectively, graphs $\G(A)$ which are complete graphs together with horns) for all bounded semirings $A$, and we also characterized all rings $R$ whose annihilating graphs $\mathbb{AG}(R)$
are either bipartite or complete together with horns.

\vs{3mm}\nin {\bf Corollary 3.5.}  {\it For any commutative ring $R$, the following statements are equivalent:}

(1) {\it Either $R\cong F_1\times F_2$ for some fields $F_i$, or $R$ is local with two nontrivial ideals.}

(2) {\it The annihilating ideal graph $\mathbb{AG}(R)$ is isomorphic to the complete graph $K_2$.}

(3) {\it Either $R\cong F_1\times F_2$ for some fields $F_i$, or $R$ is a local ring with the maximal ideal $J(R)=R\al$ for some $\al$ satisfying $\al^3=0$ and $\al^2\not=0$.}

\Proof $(1)\Lra (2):$  Clear.

$(2)\Lra (3):$ Assume $\mathbb{AG}(R)\cong K_2$. Then $R$ has exactly two proper ideals by \cite[Theorem 1.4]{BR}. In particular, $R$ is artinian and thus
the ideal $J(R)$ is finitely generated. It follows that $J(R)=R\al$ for some $\al\in J(R)$.

If $J(R)^2=0$, then either $R\cong F_1\times F_2$ for some fields $F_i$, or $R$ is a local ring by Lemma 3.1(2) and Corollary 3.4. If further $R$ is a local ring, then $J(R)=U(R)\al$ since $J(R)\al=0$. Under the assumption, $R$ has exactly one nontrivial ideal, a contradiction. The contradiction shows $R\cong F_1\times F_2$ if $J(R)^2=0$.

Now assume $J(R)^2\not=0$. By Corollary 3.4,  $R$ is a local ring with exactly two nontrivial ideals. Since $J(R)^2\not=J(R)$ by Nakayama Lemma,  $J(R)$ and $J(R)^2$ are the nontrivial ideals of $R$ by Theorem 3.3(2). Hence $J(R)^3=0$ and $J(R)^2\not=0$. In the case, $R$ is a local ring whose unique maximal ideal is $J(R)=R\al$, where $\al^3=0, \al^2\not=0$.

$(3)\Lra (1):$  Now assume that $R$ is a local ring with $J(R)=R\al$, where $\al^3=0,\al^2\not=0$. Clearly, $R\al\not=R\al^2$ and $R\al\cdot R\al^2=0$. $R=U(R)\cup R\al$ and hence $R\al=U(R)\al\cup U(R)\al^2\cup\{0\}$. Thus for any $\beta\in J(R)$, $R\beta=R\al$ if $\beta\in U(R)\al$, while $R\beta=R\al^2$ if $0\not=\beta\in U(R)\al^2$. It shows that $R\al$ and $R\al^2$ are all the possible nontrivial ideals of $R$. \quad\qed

\vs{3mm}The complete isomorphic classification of finite local rings with at most three nontrivial ideals will be discussed in a separate paper, see section four of \cite{ WL2}.

\vs{3mm}
 Note that a bounded semiring $A$ with $|Z(A)|=1$ always satisfies condition $(C_3)$. By Theorems 3.2 and 3.3, {\it the structure of a bounded semiring $A$ is completely determined by structures of integral bounded semirings if  $A$ satisfies one of the following conditions:}

(1) $|Z(A)|=1$.

(2) $|Z(A)|=2$ and $Z(A)^2\not=0$.

 On the other hand, the structure of $A$ with $Z(A)^2=0$ seems to be a little bit more complicated. Let $A$ be a bounded semiring with $Z(A)=\{c,u\}$ and $Z(A)^2=0$. Set $A_1=A\setminus Z(A)$. Then $A_1$ is an integral sub-bounded semiring (see Proposition 3.6 below). Clearly,  $u<x+u<1$ for any $x\not\in \{0,1,c,u\}$. By Lemma 3.1(2.2), $c$ is  the least nonzero element of $A$, so $cA_1^*=\{c\},c+x=x,\,\forall x\in A^*$. But it is hard to determine the partial order between $u$ and elements of $A_1^*$. By Examples 3.6 and 3.7 below, it seems that Lemma 3.1(2.2) is the best possible result.

\vs{3mm}\nin{\bf Proposition 3.6.}  {\it $A$ is a bounded semiring with condition $(C_3)$, $|Z(A)|=2$ and $Z(A)^2=0$  if and only if there exist  an integral bounded semiring $A_1$ and elements $c,u\not\in A_1$, such that $A\cong Z(A)\cup A_1$, where $Z(A)=\{c,u\}$, the addition and multiplication in $Z(A)\cup A_1$ are extended respectively from that of the bounded semiring $A_1$, and the following conditions are fulfilled:}

(1) {\it $0+x=x \,(\forall x\in Z(A)),c+y=y \,(\forall y\in A^*), u+1=1,u+u=u$, and $u+(x+y)=(u+x)+y,\,u+(u+x)=u+x\,(\forall x,y\in A_1^*)$.}

(3) {\it $Z(A)^2=0, 0Z(A)=0, cx=c\,(\forall x\in A_1^*),0\not\in u\cdot A_1^*$.}

(4) {\it $xu=u$ and $yu=u$ implies $(xy)u=u$.}

(5) {\it $x\ge y$ in $A_1$ and $yu=u$ implies $xu=u$.}

(6) {\it For any $y,z\in A_1$, $x(y+u)=xy+xu$. Also, $uy=c=uz$ implies $u(y+z)=c$.}

\Proof $\Lra$: By Lemma 3.1, it suffices to show that $A_1=A\setminus Z(A)$ is a sub-bounded semiring. In fact, if $xy\in \{c,u\}$ for some $x,y\in A_1^*$,
then we have $xc=x(yc)=(xy)c=0$, giving a contradiction.

$\Lla$: We omit the detailed verification here.\quad\qed

\vs{3mm} Next  we use Proposition 3.6 to construct some bounded semirings $A$ with  $|Z(A)|=2,\,Z(A)^2=0$.

\vs{3mm}\nin{\bf Example 3.7.} Let $A=\{0,1,c,u,b_1,b_2,\cdots\}$ be a poset with $4\le |A|\le \aleph_0$, where $0<c<u<b_1<b_2<\cdots$. Define
the binary addition by max operation. Define three commutative multiplications by
$$0\cdot x=0,\,1\cdot x=x\,(\forall x\in A),\,cu=0=c^2,\,u^2\in\{0,c,u\}, \,y\cdot z=min\{y,z\}(\text{for other $y,z$}).$$
Then it is easy to check that $A$ is a bounded semiring for each multiplication. Clearly, $|Z(A)|=2,Z(A)^2=0$.

\vs{3mm}\nin{\bf Example 3.8.} Let $Z(A)=\{c,u\}$ and $A=Z(A)\cup A_1$, where $A_1=\{0,1,b_1,b_2,\cdots\}$ is a chain  $0<b_1<b_2<\cdots<1$, with $max$ as the binary addition and with multiplication $b_ib_j=b_1$ ($\forall i,j\ge 1 $). Then clearly $A_1$ is an integral bounded semiring. Fix a positive integer $n>1$ and extend the partial order of $A_1$ to $A$ by $0<c<u<b_n$. Extend the commutative addition to $A$ by
$$0+x=x,\,1+x=1 \, , x+x=x\, (\forall x\in \{c,u\}, c+y=y\,(\forall y\in A^*)$$ and $$u+b_i=b_n \,(1\le i\le n-1),u+b_j=b_j\,(j\ge n).$$
Extend the commutative multiplication to $A$ by
$$0Z(A)=0=Z(A)^2,\, ub_i=c\,(\forall i).$$
Note that $ c+y=y\,(\forall y\in A^*)$ implies that $c$ is the least element of $A$, thus $c^2=0$ means that condition $(C_3)$ holds. Note that $0+x=x,\,x+1=1$ implies $0\le x\le 1$. Then it is easy to apply Proposition 3.6  check that $A$ is a bounded semiring with $|Z(A)|=2, Z(A)^2=0$. {\it Note that $u$ and $b_i$ are incomparable for any $i$ with $1\le i\le n-1$, while by Lemma 3.1(2), $u< p$ for any prime element $p$ with $p\not=c,u$.}

\vs{3mm}For any bounded semiring $A$, recall that $A^{(r)}$ is the direct product of $r$ copies of $A$. We have the following.

\vs{3mm}\nin{\bf Proposition 3.9.}  {\it  $A$ is a bounded semiring, in which both DCC and condition $(C_3)$ holds for elements of $A$, if and only if there exist an integer $n\ge 1$ and a bounded semiring $A_1$ such that either $A\cong A_1$ or $A\cong \{0,1\}^{(n)}\times A_1$, where $c^2=0$ for each minimal element $c$ of $A_1$ and DCC holds for elements of $A_1$.  }

\Proof $\Lla:$ Let $A_1$ be a bounded semiring, and assume that DCC holds for elements of $A_1$ and $c^2=0$ for each minimal element $c$ of $A_1$. Then there exists at leat one minimal element $c$ in $A$ and hence $Z(A_1)\not=\eset$. If $A\cong A_1$, then clearly condition $(C_3)$ holds in $A$.  If $A\cong \{0,1\}^n\times A_1$ for some finite $n\ge 1$,  then both DCC and condition $(C_3)$ also hold for elements of $A$.

$\Lra:$ Assume that DCC holds for elements of $A$. Then there exists at leat one minimal element in $A$.  If further $c^2=0$ for each minimal element $c$ of $A$, then $Z(A)\not=\eset$ and condition $(C_3)$ holds in $A$. Then we take $A_1=A$. In the following assume that there exists an idempotent minimal element $e_1$ in $A$. Then by condition $(C_3)$, there exists a nonzero idempotent $f_1\in A$ such that $A\cong \{0,1\}\times Af_1$.
Clearly,  both DCC and condition $(C_3)$ hold for elements of the bounded semiring $Af_1$ and an induction shows that $A\cong \{0,1\}^{(n)}\times A_1$ for some $n\ge 1$, where $A_1$ is a bounded semiring in which $c^2=0$ for each minimal element $c$ of $A_1$ and DCC holds for elements of $A_1$.
\quad\qed


\begin{thebibliography}{gg}

\bibitem{AA}G. Aalipour, S. Akbari, R. Nikandish,
M.J. Nikmehr and F. Shaveisi. Minimal prime ideals and cycles in
annihilating-ideal graphs, {\bf Rocky Mountain J. Math.}   {\bf 43:5}(2013) $1415-1425$.

\bibitem{AA2} G. Aalipour, S. Akbari, M. Behboodi, R. Nikandish,
M.J. Nikmehr and F. Shaveisi. The classification of the annihilating-ideal graph
of a commutative ring. {\bf Algebra Colloq.} (accepted)

\bibitem{AL}
D.F. Anderson and P.S. Livingston.  The zero-divisor graph of a
commutative ring, {\bf  J. Algebra} {\bf 217}(1999) $434-447$.

\bibitem{AM}
M.F. Atiyah and I.G. MacDonald. {\it
Introduction to Commutative Algebra}, {\bf Addison-Wesley}, Reading, MA, $1969$.

\bibitem{BR} M. Behboodi and Z. Rakeei. The annihilating-ideal graph of commutative
rings I, {\bf Journal of
Algebra and its Applications} {\bf 10:4}(2011) $727-739$.

\bibitem{BMR} S. Bistarelli, U. Montanari,
and F. Rossi. Semiring-based constraint satisfaction and
optimization, {\bf Journal of the ACM} {\bf 44:2}(1997) $201-236$.

\bibitem{BMRSVF}  S. Bistarelli,  U. Montanari,  F. Rossi, Schiex, T.;  G.  Verfaillie,  H. Fargier. Semiring-based CSPs and valued CSPs: frameworks, properties, and comparison, {\bf Constraints} {\bf 4:3} (1999) $199-240$.
    
\bibitem{Bjorner} A. Bj\"orner, Shellable and Cohen-Macaulay partially ordered sets. {\bf Trans. Amer. Math.}
Soc., Vol.260 (1980), pp.159-183.

\bibitem{DMS}
F.R. DeMeyer, T. McKenzie, and K. Schneider. The zero-divisor graph
of a commutative semigroup, {\bf Semigroup Forum} {\bf 65}(2002)
$206-214$.

\bibitem{Lam}T.Y. Lam. {\it A First Course in Non-Commutative Rings}, {\bf Springer-Verlag}, New York, Inc
$1991$.

\bibitem{LW}D.C. Lu and T.S. Wu. On the normal ideals of exchange rings, {\bf  Siberian Math.J.} {\bf 49}(2008) $663-668$.

\bibitem{McL}K.R. McLean. Commutative artinian principal ideal rings. {\bf Proc. London Math. Soc.} {\bf 26:3} (1973), $249-272$.

\bibitem{WL}  T.S. Wu and  D.C. Lu. Sub-semigroups determined by the zero-divisor graph,
{\bf  Discrete Math.} {\bf 308}(2008) $5122-5135$.

\bibitem{WL2}  T.S. Wu, H.Y. Yu and  D.C. Lu. The structure of finite local principal ideal rings. {\bf Comm. Algebra} {\bf 40:12}(2012) 4727-4738.

\bibitem{YW}  H.Y. Yu and T.S. Wu. On realizing zero-divisor graphs of bounded semirings. Preprint.

\end{thebibliography}
\end{document}